%version of August 19 , 2008
%First draft: completed August 22, 2008. 
% Second draft: August 30, 2008.
%Final draft of September 8, 2008
\input amstex

\documentstyle{amsppt} 
\magnification1200
\NoBlackBoxes
\pagewidth{6.5 true in}
\pageheight{9.25 true in} 
\document

\topmatter 
\title
Mass  equidistribution for Hecke eigenforms 
\endtitle
\author 
R. Holowinsky and K. Soundararajan
\endauthor
\address 
Department of Mathematics, University of Toronto, 
Toronto, Ontario, Canada M5S 2E4
\endaddress 
\email 
romanh{\@}math.toronto.edu
\endemail
\address
{Department of Mathematics, 450 Serra Mall, Bldg. 380, Stanford University, 
Stanford, CA 94305-2125, USA}
\endaddress
\email 
ksound{\@}stanford.edu
\endemail
\thanks   The first author is supported by NSERC and the Fields Institute, Toronto.  The second author is partially supported by the National Science Foundation (DMS 0500711) 
and the American Institute of Mathematics (AIM). 
\endthanks
\endtopmatter

\def\lam{\lambda}
\def\Lam{\Lambda}
\document  

\head 1.  Introduction \endhead

\noindent A central problem in the area of quantum chaos is 
to understand the limiting behavior of eigenfunctions.  An important example which 
has attracted a lot of attention is the case of Maass cusp forms of large Laplace eigenvalue 
for the space $X=SL_2({\Bbb Z})\backslash {\Bbb H}$.  Let $\phi$ denote such a 
Maass form, and let $\lam$ denote its Laplace eigenvalue, and let $\phi$ be 
normalized so that $\int_X |\phi(z)|^2 \frac{dx \ dy}{y^2}=1$.  Zelditch [19] has 
shown\footnote{We have given Zelditch's result in the context of $SL_2({\Bbb Z})\backslash {\Bbb H}$.
In fact he proves more, since he considers equidistribution of the micro-local 
lift to $SL_2({\Bbb Z})\backslash SL_2({\Bbb R})$.  Lindenstrauss's result stated below 
also holds for this micro-local lift.} that as $\lam\to \infty$, for a {\sl typical} Maass form $\phi$ the measure $\mu_{\phi}:= |\phi(z)|^2 \frac{dx \ dy}{y^2} $ approaches 
the uniform distribution measure $\frac 3{\pi} \frac{dx \ dy}{y^2}$. 
% For example,  
%consider the space $X=SL_2({\Bbb Z}) \backslash {\Bbb H}$ 
%and let $\phi$ be a Maass cusp form with Laplace eigenvalue $\lambda$.  
This statement is referred to as ``Quantum Ergodicity."   Rudnick and Sarnak [13] 
have conjectured that an even stronger result holds.  Namely, that as $\lam \to \infty$, 
for {\sl every} Maass form $\phi$ the measure $\mu_{\phi}$ approaches the 
uniform distribution measure.  This conjecture 
is referred to as ``Quantum Unique Ergodicity."  
%Suppose also that $\phi$ is an eigenfunction of all the Hecke operators, and that 
%$\phi$ has been normalized so that $\int_{X} |\phi(z)|^2 \frac{dx \ dy}{y^2}=1$.   
%Rudnick and Sarnak [9] have then conjectured (quantum unique ergodicity) 
%that for large $\lambda$ 
%the measure $\mu_{\phi}:= |\phi(z)|^2 \frac{dx \ dy}{y^2} $ approaches 
%the uniform distribution measure $\frac 3{\pi} \frac{dx \ dy}{y^2}$.  
Lindenstrauss [8] has made great progress towards this conjecture, 
showing that, for Maass cusp forms that are eigenfunctions of the Laplacian and 
all the Hecke operators,\footnote{The spectrum of the Laplacian is expected to be simple, 
so that any eigenfunction of the Laplacian would automatically be an eigenfunction of 
all Hecke operators.  This is far from being proved.}  the only possible limiting measures are of the form 
$\frac 3{\pi} c \frac{dx \ dy}{y^2}$ with $0\le c\le 1$.   For illuminating 
accounts on this conjecture we refer the reader to [7, 8, 9, 10, 11, 13, 14, 15, 18].

Here we consider a holomorphic analog of the   
 quantum unique ergodicity conjecture.  This analog is very much in the 
 spirit of the Rudnick-Sarnak conjectures, and has been spelt out explicitly in [10, 14].  
 Let $f$ be a holomorphic modular cusp form of weight $k$ (an even integer) for $SL_2({\Bbb Z})$.  
 Associated to $f$ we have the %, and suppose that 
% $f$ is an eigenfunction of all the Hecke operators.   The invariant 
 measure %here is 
 $$
 \mu_f := y^k |f(z)|^2 \frac{dx \ dy}{y^2},
 $$ 
 which is invariant under the action of $SL_2({\Bbb Z})$, and 
 we suppose that $f$ has been normalized so that 
 $$ 
 \int_{X} y^{k} |f(z)|^2 \frac{dx \ dy}{y^2} = 1.
 $$ 
 The space $S_{k}(SL_2({\Bbb Z}))$ of cusp forms of weight $k$ for 
 $SL_2({\Bbb Z})$ is a vector space of dimension about $k/12$, and contains 
 elements such as $\Delta(z)^{k/12}$ (if $12|k$, and where $\Delta$ is 
 Ramanujan's cusp form) for which the measure will not tend to uniform distribution.  
% cannot expect $\mu_f$ to tend to uniform distribution for all forms in this 
 %space.\footnote{ In contrast the space of Maass forms with a given eigenvalue 
 %is expected to have dimension at most one.}   
 Therefore one restricts attention 
 to a particularly nice set of cusp forms, namely those that are eigenfunctions 
 of all the Hecke operators.  The Rudnick-Sarnak conjecture in this context 
 states that as $k\to \infty$, for {\sl every} Hecke eigencuspform $f$ the 
 measure $\mu_f$ tends to the uniform distribution measure.  For 
 simplicity, we have restricted ourselves to the full modular group, but 
 the conjecture could be formulated just as well for holomorphic newforms of level $N$.  
 Luo and Sarnak [10] have shown that equidistribution holds for most Hecke eigenforms, 
 and Sarnak [14] has shown that it holds in the special case of dihedral forms.  It does not 
 seem clear how to extend Lindenstrauss's work to the holomorphic setting.\footnote{The difficulty 
 from the ergodic point of view concerns the invariance under the 
 geodesic flow of the quantum limits of the micro-local lifts associated to holomorphic forms.}

 In this paper we shall establish the Rudnick-Sarnak conjecture 
 for holomorphic Hecke eigencuspforms.  The 
 proof combines two different approaches to the mass equidistribution conjecture, 
 developed independently by the authors [4, 17].  
 Either of these approaches is capable of showing that there are very 
 few possible exceptions to the conjecture, and under reasonable 
 hypotheses either approach would show that there are no exceptions. 
 However, it seems difficult to show unconditionally that 
 there are no exceptions using just one of these approaches.  
 Fortunately, as we shall explain below, the two approaches 
 are complementary, and the few rare cases that are untreated by 
 one method fall easily to the other method.  Both approaches 
 use in an essential way that the Hecke eigenvalues of a holomorphic 
 eigencuspform satisfy the Ramanujan conjecture (Deligne's theorem).  
 The Ramanujan conjecture remains open for Maass forms, and this is 
 the (only) barrier to using our methods in the non-holomorphic 
 setting.\footnote{Assuming the Ramanujan conjecture for Maass forms, 
 our methods would obtain the stronger micro-local version of QUE.}
 At present, it is not clear how to use our methods in the case 
 of compact quotients.  
 
 Recall that for two smooth bounded functions $g_1$ and $g_2$ on $X$ 
 we may define the Petersson inner product 
 $$ 
 \langle g_1, g_2 \rangle
 = \int_{X} g_1(z) \overline{g_2(z)} \frac{dx \ dy}{y^2}.  
 $$
 In this definition we could allow for one of the $g_1$ or $g_2$ to be unbounded, 
 so long as the other function decays appropriately for the integral to converge. 
 If $f$ is a modular form 
 of weight $k$, below we shall let $F_k(z)$ denote $y^{k/2} f(z)$ where $z=x+iy$. 
  Let $h$ denote a smooth bounded function on $X$. Considering $h$ as fixed, 
  and letting $k\to \infty$, the Rudnick-Sarnak conjecture asserts that for every 
  Hecke eigencuspform $f$ of weight $k$ we have\footnote{We are slightly abusing 
  notation here, because $F_k$ is not $SL_2({\Bbb Z})$-invariant. However $|F_k(z)|^2$ 
  is $SL_2({\Bbb Z})$-invariant, and so the inner product in (1.1) does not 
  depend on a choice of the fundamental domain.}
  $$ 
  \langle h F_k, F_k \rangle \to \frac 3\pi \langle h,1\rangle, \tag{1.1}
  $$ 
with the rate of convergence above depending on the function $h$.

To attack the conjecture (1.1), it is convenient to decompose 
the function $h$ in terms of a basis of smooth functions on $X$. 
There are two natural ways of doing this.  First we could use the 
spectral decomposition of a smooth function on $X$ in terms 
of eigenfunctions of the Laplacian.  The spectral expansion 
will involve (i) the constant function $\sqrt{3/\pi}$, (ii) 
Maass cusp forms $\phi$ that are also eigenfunctions of 
all the Hecke operators, and (iii) Eisenstein series on the $\tfrac 12$ line.  
Recall that the Eisenstein series is defined for Re$(s) >1$ by 
$$ 
E(z,s) = \sum_{\gamma \in \Gamma_{\infty} \backslash \Gamma} \text{Im}(\gamma z)^s,
$$ 
where $\Gamma =SL_2({\Bbb Z})$ and $\Gamma_\infty$ denotes the 
stabilizer group of the cusp at infinity (namely the set of all 
translations by integers).  The Eisenstein series $E(z,s)$ admits a
meromorphic continuation, with a simple pole at $s=1$, and is analytic 
for $s$ on the line Re$(s)=\tfrac 12$.  For more on the spectral expansion 
see Iwaniec [5].  Note that (1.1) is trivial 
when $h$ is the constant eigenfunction.  To establish (1.1) using the spectral decomposition, 
we would need to show that for a fixed Maass eigencuspform $\phi$, and 
for a fixed real number $t$ that 
$$ 
\langle \phi F_k, F_k \rangle, \qquad \text{and} \qquad \langle   E(\cdot, \tfrac 12+it) F_k, F_k\rangle 
\to 0,
$$ 
as $k\to \infty$.  The inner products above may be related to 
values of $L$-functions.  In the case of Eisenstein series, this is 
the classical work of Rankin and Selberg.  In the more difficult 
Maass form case, this relation (to a triple product $L$-function) is given by a beautiful formula of 
Watson [18].  The connection to $L$-functions, and estimating such values, forms the 
basis for Soundararajan's approach to (1.1).  

Alternatively, one could expand the function $h$ in terms of incomplete 
Poincare and Eisenstein series.  Let $\psi$ denote a smooth 
function, compactly supported in $(0,\infty)$.  For an integer $m$ the 
incomplete Poincare series is defined by 
$$ 
P_m(z \mid \psi) = \sum_{\gamma \in \Gamma_\infty \backslash \Gamma} e(m\gamma z)\psi( \text{Im}(\gamma z)).
$$ 
In the special case $m=0$ we obtain incomplete Eisenstein series $E(z\mid \psi) = P_0(z\mid \psi)$.  
For an account on approximating a smooth function $h$ using incomplete Poincare series see [9]. 
Luo and Sarnak [10] noted that in this approach to (1.1) one faces the problem of estimating the 
shifted convolution sums (for $m$ fixed, and as $k\to \infty$) 
$$
\sum_{n \asymp k} \lam_f(n) \lam_f(n+m),
$$ 
where the sum is over $n$ of size $k$, and $\lam_f(n)$ denotes the Hecke eigenvalue 
of $f$ normalized so that Deligne's bound reads $|\lam_f(n)| \le d(n)$.  The study 
of these shifted sums using sieve methods forms the basis for Holowinsky's 
approach to (1.1).

We are now in a position to state our main result, after which we 
will describe the main Theorems of Holowinsky and Soundararajan, 
and how those combine.  
 
 \proclaim{Theorem 1}  Let $f$ be a Hecke eigencuspform of weight $k$ for $SL_2({\Bbb Z})$, 
 and write $F_k(z)=y^{k/2} f(z)$.  
 
 \noindent (i).  Let $\phi$ be a Maass cusp form which is also an eigenfunction 
 of all Hecke operators.  Then 
 $$ 
| \langle \phi F_k, F_k  \rangle  | \ll_{\phi, \epsilon} (\log k)^{-\frac {1}{30}+\epsilon}.
$$ 

\noindent (ii).   Let $\psi$ be a fixed smooth function compactly supported in $(0,\infty)$.  Then 
$$ 
\Big| \langle   E(\cdot \mid \psi) F_k, F_k\rangle - \frac 3{\pi} \langle  E(\cdot \mid \psi), 1\rangle 
\Big| \ll_{\psi, \epsilon} (\log k)^{-\frac 2{15}+\epsilon}.
$$ 
\endproclaim 

\demo{Remark 1}
We have made no attempt to optimize the rate of decay given above.  Our 
methods would not appear to lead to any decay better than $(\log k)^{-1}$.  
The generalized Riemann hypothesis is known to imply a rate of decay $k^{-\frac 12+\epsilon}$ 
which would be optimal, see [7, 11].  
\enddemo

\demo{Remark 2} 
A striking consequence of Theorem 1 is that the zeros of the modular 
form $f$ lying in $X$ (there are about $k/12$ such zeros) become equidistributed 
with respect to the measure $\frac 3{\pi} \frac{dx \ dy}{y^2}$.  This follows from the 
work of Rudnick [12], who derived this consequence from the mass  equidistribution conjecture.  

\enddemo

\demo{Remark 3}  With more effort, we could keep track of the dependence on $\phi$ and $\psi$ 
in Theorem 1.  This would require some careful book-keeping in the works [4] and [17].   Keeping 
track of these dependencies would allow one to give a rate of decay for the discrepancy (for example, the spherical cap discrepancy defined in [9]) between 
the measure $\mu_f$ and the uniform distribution measure.  

\enddemo

%Our proof of  Theorem 1 combines two different approaches of the 
%authors developed independently.  Both approaches invoke the Ramanujan 
%conjectures in an essential way.  In the holomorphic case this is known due to Deligne, but 
%these bounds are unavailable in the case of Maass forms.  Below, $\lam_f(n)$ denotes the 
%eigenvalue of the $n$-th Hecke operator 
% $T_n$ acting on $f$.  We normalize the Hecke operators so that Deligne's bound 
% reads $|\lam_f(p)| \le 2$.  
 
 We now describe the results of Holowinsky and Soundararajan, and how they pertain 
 to Theorem 1.  In [4] Holowinsky attacks the inner products in Theorem 1 by an unfolding method, 
 which leads him to the estimation of shifted convolution sums of the Hecke eigenvalues. 
He then develops a sieve method to estimate these shifted convolution sums,   
obtaining the following result.  

\proclaim{Theorem 2 (Holowinsky)}  Keep the notations of Theorem 1.  Let $\lam_f(n)$ 
denote the Hecke eigenvalue of $f$ for the $n$-th Hecke operator normalized so that 
$|\lam_f(n)|\le d(n)$.  Let $L(s,\text{\rm sym}^2 f)$ denote the symmetric square $L$-function 
attached to $f$, and define  
$$ 
M_k(f):= \frac{1}{(\log k)^2L(1,\text{\rm sym}^2 f)} \prod_{p\le k} \Big(1+\frac{2|\lam_f(p)|}{p}\Big).
$$  

\noindent (i). For a Maass cusp form $\phi$ we have 
$$ 
| \langle \phi F_k, F_k  \rangle | \ll_{\phi,\epsilon} (\log k)^{\epsilon} M_k(f)^{\frac 12}.
$$ 

\noindent (ii). For an incomplete Eisenstein series $E(z\mid \psi)$ we have 
$$
\Big| \langle   E(\cdot \mid \psi)F_k, F_k \rangle - \frac 3{\pi} \langle E(\cdot \mid \psi), 1\rangle 
\Big| \ll_{\psi,\epsilon} (\log k)^{\epsilon} M_k(f)^{\frac 12} (1+R_k(f)), 
$$
where
$$ 
R_k(f) = \frac{1}{k^{\frac 12}L(1,\text{\rm sym}^2 f)} \int_{-\infty}^{+\infty} \frac{|L(\tfrac 12+it, \text{\rm sym}^2 f)|}{(1+|t|)^{10}} dt.
$$
\endproclaim 

Although this is not immediately apparent, the quantity $M_k(f)$ appearing above is 
expected to be small in size.   One can show that there are at most $K^{\epsilon}$ 
eigenforms $f$ with weight below $K$ for which $M_k(f) \ge (\log k)^{-\delta}$ 
for some fixed $\delta >0$.  A weak form of the generalized Riemann hypothesis 
could be used to show that $M_k(f) \le (\log k)^{-\delta}$ for some $\delta>0$.  
Moreover one can show that 
$$ 
M_k(f)\ll (\log k)^{\epsilon} \exp\Big( -\sum_{p\le k} \frac{(|\lam_f(p)|-1)^2}{p}\Big), \tag{1.2}
$$
so that one would expect $M_k(f)$ to be small unless $|\lam_f(p)| \approx 1$ for 
most $p\le k$.  This last possibility is not expected to hold,\footnote{Since we would 
expect $\lam_f(p)$ to be distributed in $[-2,2]$ according to the Sato-Tate 
measure.}  and can be shown to be 
rare,\footnote{By bounding the quantity in (1.2) in terms of $L(1,\text{sym}^2 f)$ and $L(1,\text{sym}^4 f)$; 
see [4] for details.}  but is difficult to rule out completely.   In the Eisenstein series case ((ii) above), 
one also needs to bound $R_k(f)$; again this can be shown to be small in 
all but very rare cases.

\smallskip 
As mentioned earlier, if we approach (1.1) through the spectral expansion of $h$ we 
are led to estimating central values of $L$-functions.  Here it turns out that an 
easy convexity bound for $L$-values barely fails to be of use, and improved 
subconvexity estimates (saving a power of the analytic conductor) would solve the problem completely
(see for example [7]).  
In [17] Soundararajan developed a general method which gives weak subconvexity 
bounds for central values of $L$-functions.  Instead of obtaining a power 
saving of the analytic conductor, one obtains a saving of powers of the 
logarithm of the analytic conductor.

  \proclaim{Theorem 3 (Soundararajan)}  Keep the notations of Theorems 1 and 2.
  
  \noindent (i).  For a Maass eigencuspform  $\phi$ we have 
  $$ 
|  \langle \phi F_k, F_k  \rangle | \ll_{\phi,\epsilon} \frac{(\log k)^{-\frac 12+\epsilon} }{L(1,\text{\rm sym}^2 f)}.
$$ 

\noindent (ii).  For the Eisenstein series $E(z,\tfrac 12+it)$ we have 
$$
|\langle  E(\cdot, \tfrac 12+it) F_k, F_k \rangle | \ll_{\epsilon} (1+|t|)^2 \frac{(\log k)^{-1+\epsilon}}{L(1,
\text{\rm sym}^2 f)}.
$$
\endproclaim

If $L(1,\text{sym}^2 f) \ge (\log k)^{-\frac 12+\delta}$ for some $\delta>0$ then 
Theorem 3 would establish the Rudnick-Sarnak conjecture.  This bound on $L(1,\text{sym}^2 f)$ 
is certainly expected to hold; for example it follows from a weak form of the generalized Riemann hypothesis.  Moreover, one can show that this bound fails to hold for at most $K^{\epsilon}$ 
eigencuspforms $f$ with weight below $K$.   However, 
we know only that $L(1,\text{sym}^2 f) \gg (\log k)^{-1}$, and it seems difficult to rule out the 
possibility of small values of $L(1,\text{sym}^2 f)$ completely.  

\smallskip

We have seen how either of the approaches in Theorems 2 and 3 should work 
always, and that both fail in rare circumstances which are difficult to rule 
out.   Now we shall see that in the rare circumstances that one 
of these results fails, the other result succeeds.  To gain an intuitive understanding 
of this phenomenon, note that Theorem 3 fails only when $L(1,\text{sym}^2 f) \le 
(\log k)^{-\frac 12+\delta}$ is small.  But this $L$-value is small only if for most 
primes $p\le k$ we have $\lam_f(p^2) \approx -1$ (a Siegel zero type phenomenon).  
But this means $\lam_f(p)^2 -1 \approx -1$ so that $\lam_f(p)\approx 0$.   Recall 
now that the quantity $M_k(f)$ appearing in Holowinsky's work is small unless 
$|\lam_f(p)| \approx 1$ for most $p\le k$.  Evidently, both situations cannot happen 
simultaneously. 
   
More precisely, in Lemma 3 below we shall show that 
$$ 
M_k(f) \ll (\log k)^{\frac 16}  (\log \log k)^{\frac 92} L(1,\text{sym}^2 f)^{\frac 12}.
$$
Therefore, if $L(1,\text{sym}^2 f) < (\log k)^{-\frac 13 -\delta}$ for 
some small $\delta >0$ it follows from Theorem 2 (i) that $\langle \phi F_k, F_k \rangle$ is 
small.  However, if $L(1,\text{sym}^2 f) > (\log k)^{-\frac 13 -\delta} >(\log k)^{ -\frac 12+\delta}$ 
then as noted above Theorem 3 (i) shows that $\langle \phi F_k ,F_k \rangle$ is small.    
This shows how Theorems 2 and 3 complement each other in the cusp form case.  
In the case of Eisenstein series, we first show how the weak subconvexity 
results in [17] lead to  a satisfactory bound for the term $R_k(f)$ appearing 
in Theorem 2 (ii) (see Lemma 1 below).  Then the argument follows as in case (i).

{\bf Acknowledgments.}  We are grateful to Peter Sarnak for 
encouragement and several valuable discussions.  
  
 \head 2.  Preliminaries on the symmetric square $L$-function \endhead
 
\noindent We collect together here some properties of $L(s,\text{sym}^2 f)$ that 
we shall need for the proof of Theorem 1.   We shall require several important 
results on $L(s,\text{sym}^2 f)$ that are known thanks to the works of Shimura [16], 
Gelbart and Jacquet [1],  Hoffstein and Lockhart [3], and Goldfeld, Hoffstein and Lieman [2].  
If we write 
$$ 
L(s,f) =\sum_{n=1}^{\infty} \frac{\lam_f(n)}{n^s} = \prod_p 
\Big( 1-\frac{\alpha_p}{p^s}\Big)^{-1} \Big( 1-\frac{\beta_p}{p^s}\Big)^{-1}, 
$$ 
where $\alpha_p$ and $\beta_p=\overline{\alpha_p}$ are complex numbers of magnitude $1$, 
then the symmetric square 
$L$-function is 
$$
L(s,\text{sym}^2 f) = \sum_{n=1}^{\infty} \frac{\lam^{(2)}_f(n)}{n^s} = 
\prod_{p} \Big(1-\frac{\alpha_p^2}{p^s}\Big)^{-1} \Big(1-\frac 1{p^s}\Big)^{-1} 
\Big(1-\frac{\beta_p^2}{p^s}\Big)^{-1}. 
$$ 
The series and product above converge absolutely in Re$(s)>1$, and 
by the work of Shimura [16], we know that $L(s, \text{sym}^2 f)$ extends 
analytically to the entire complex plane, and  satisfies the functional equation 
$$
\Lam(s,\text{sym}^2 f) = \Gamma_{\Bbb R}(s+1)\Gamma_{\Bbb R}(s+k-1) 
\Gamma_{\Bbb R}(s+k) L(s,\text{sym}^2 f) = \Lam(1-s,\text{sym}^2 f), \tag{2.1}
$$
where $\Gamma_{\Bbb R}(s) = \pi^{-s/2} \Gamma(s/2)$.  

Gelbart and Jacquet [1] have shown that $L(s,\text{sym}^2 f)$ arises as the 
$L$-function of a cuspidal automorphic representation of $GL(3)$.  Therefore, 
invoking the Rankin-Selberg convolution for $\text{sym}^2 f$, one can 
establish a classical zero-free region for $L(s,\text{sym}^2 f)$.  For 
example, from Theorem 5.42 (or Theorem 5.44) of Iwaniec and Kowalski [6] one obtains 
that for some constant $c>0$ the region 
$$
{\Cal R} = \Big\{ s=\sigma+it: \ \ \sigma \ge 1-\frac{c}{\log k(1+|t|)} \Big\} 
$$
does not contain any zeros of $L(s,\text{sym}^2 f)$ except possibly for a 
simple real zero.   The work of Hoffstein and Lockhart [3] (see the 
appendix by Goldfeld, Hoffstein and Lieman [2]) shows that $c>0$ may 
be chosen so that there is no real zero in our region ${\Cal R}$.  Thus 
$L(s,\text{sym}^2 f)$ has no zeros in ${\Cal R}$.   Moreover the 
work of Goldfeld, Hoffstein, and Lieman [2] shows that 
$$ 
L(1,\text{sym}^2 f )\gg \frac{1}{\log k}. \tag{2.2} 
$$ 
To be precise, the work of Goldfeld, Hoffstein, and Lieman considers 
symmetric square $L$-functions of Maass forms in the eigenvalue 
aspect, but our case is entirely analogous, and follows upon making 
minor modifications to their argument. 

\proclaim{Lemma 1}  For any $t\in {\Bbb R}$ we have 
$$ 
|L(\tfrac 12+it,\text{\rm sym}^2 f)| \ll \frac{k^{\frac 12} (1+|t|)^{\frac 34}}{(\log k)^{1-\epsilon}}.
$$
Therefore the quantity $R_k(f)$ appearing in Theorem 2 (ii) satisfies 
$$ 
R_k(f) \ll \frac{(\log k)^{\epsilon}}{(\log k) L(1,\text{\rm sym}^2 f)} \ll (\log k)^{\epsilon}.
$$ 
\endproclaim 

\demo{Proof}  The bound on $L(\tfrac 12+it, \text{sym}^2 f)$ 
 follows from the results of [17] on weak subconvexity; see Example 1 
there.  Using this bound in the definition of $R_k(f)$ immediately 
gives the stated estimate.  
\enddemo

\proclaim{Lemma 2}  We have 
$$ 
L(1,\text{\rm sym}^2 f) \gg (\log \log k)^{-3} \exp\Big( \sum_{p\le k} \frac{\lam_f(p^2)}{p}\Big).
$$
\endproclaim

\demo{Proof} Let $1\le \sigma \le \frac 54$, and consider for some $c>1$, and 
$x\ge 1$
$$ 
\frac{1}{2\pi i} \int_{c-i\infty}^{c+i\infty} -\frac{L^{\prime}}{L}(s+\sigma, \text{sym}^2 f) \frac{2x^s}{s(s+2)} ds. \tag{2.3}  
$$
We write 
$$
-\frac{L^{\prime}}{L}(s,\text{sym}^2 f) = \sum_{n=1}^{\infty} \frac{\Lam_{\text{sym}^2 f}(n)}{n^s}
$$
say where $\Lam_{\text{sym}^2 f}(n)=0$ unless $n=p^k$ is a prime power, in 
which case $\Lam_{\text{sym}^2 f}(p) = \lam_f(p^2)$ and for $k\ge 2$, $|\Lam_{\text{sym}^2 f}(p^k)|\le 3$.
Using this in (2.3) and integrating term by term we see that 
$$ 
\align
\frac{1}{2\pi i} \int_{c-i\infty}^{c+i\infty} -\frac{L^{\prime}}{L}(s+\sigma, \text{sym}^2 f) \frac{2x^s}{s(s+2)}
 ds
&= \sum_{p\le x } \frac{\lam_f(p^2)\log p}{p^{\sigma}} \Big( 1- \Big(\frac px\Big)^2 \Big) + O(1)\\
&= \sum_{p\le x} \frac{\lam_f(p^2) \log p}{p^{\sigma}} +O(1). \tag{2.4}\\
\endalign
$$ 

We next evaluate the LHS of (2.4) by shifting the line of integration to 
the line Re$(s)= -\frac 32$.  We encounter poles at $s=0$, and when $s= \rho-\sigma$ 
for a non-trivial zero $\rho=\beta+i\gamma$ of $L(s,\text{sym}^2 f)$.  Computing these 
residues, we obtain that the quantity in (2.4) equals 
$$
 -\frac{L^{\prime}}{L}(\sigma, \text{sym}^2 f) +O\Big(\sum_{\rho} \frac{x^{\beta-\sigma} }{|\rho-\sigma||\rho-\sigma+2|}\Big) + \frac{1}{2\pi i } \int_{-\frac 32-i\infty}^{-\frac 32+i\infty} -\frac{L^{\prime}}{L}(s+\sigma,\text{sym}^2 f) \frac{2x^s ds}{s(s+2)}. 
 $$
 To estimate the integral above, we differentiate the functional equation (2.1) 
 logarithmically, and use Stirling's formula.  Thus if $s=-\frac 32+it$ we obtain 
 that 
 $$ 
 -\frac{L^{\prime}}{L}(s+\sigma,\text{sym}^2 f) \ll \log (k(1+|t|)) + \Big| \frac{L^{\prime}}{L} (1-s-\sigma, 
 \text{sym}^2 f)\Big| \ll \log (k(1+|t|)).
 $$ 
 Thus we deduce that  
 $$ 
 \sum_{p\le x} \frac{\lam_f(p^2) \log p}{p^{\sigma}} 
 = -\frac{L^{\prime}}{L}(\sigma, \text{sym}^2 f) +O\Big(\sum_{\rho} \frac{x^{\beta-\sigma} }{|\rho-\sigma||\rho-\sigma+2|}\Big) +O(x^{-\frac 32} \log k). \tag{2.5}
 $$
 
 To bound the sum over zeros above, we split the sum into intervals where $n \le |\gamma| <n+1$ 
 for $n =0$, $1$, $\ldots$.   Each such interval contains at most $\ll \log (k(1+n))$ zeros, 
 and moreover they all lie outside the zero-free region ${\Cal R}$.  Therefore the 
 sum over zeros in (2.5) is 
 $$ 
 \ll x^{-c/\log k}  (\log k)^2 +  \sum_{n=1}^{\infty} x^{-c/\log (k(1+n))} \frac{\log (k(1+n))}{n^2} 
 \ll x^{-c/(2\log k)} (\log k)^2 + 1.
 $$ 
 Choose $x= k^{4(\log \log k)/c}$ so that the above becomes $\ll 1$.  
 Combining this estimate with (2.5) we conclude that for this choice of $x$, 
 $$ 
 -\frac{L^{\prime}}{L}(\sigma, \text{sym}^2 f) = \sum_{p\le x} \frac{\lam_f(p^2) \log p}{p^{\sigma}} + O(1),
 $$ 
 and integrating both sides from $\sigma=1$ to $\frac 54$ that 
$$
\log L(1, \text{sym}^2 f) = \sum_{p\le x} \frac{\lam_f(p^2)}{p} + O(1).
$$
Since $\sum_{k< p\le x} {\lam_f(p^2)}/p \le \sum_{k<p\le x} 3/p = 3 \log \log \log  x+ O(1)$, 
our Lemma follows.

\enddemo 

\demo{Remark}  The factor $(\log \log k)^{-3}$ above is extraneous.  With more 
effort one should be able to show that $L(1,\text{sym}^2 f) \asymp \exp(\sum_{p\le k} \lam_f(p^2)/p)$.
\enddemo

From Lemma 2 we shall obtain an estimate for the quantity $M_k(f)$ appearing in 
 Theorem 2.

\proclaim{Lemma 3}  We have 
$$ 
M_k(f) \ll (\log k)^{\frac 16} (\log \log k)^{\frac 92} L(1,\text{\rm sym}^2 f)^{\frac 12}.
$$
\endproclaim  
\demo{Proof} From the inequality $2|x| \le \frac 23  + \frac 32 x^2$, and the Hecke relations, we obtain 
that 
$$ 
 2\sum_{p\le k} \frac{|\lam_f(p)|}{p} 
\le \frac{2}{3} \sum_{p\le k} \frac 1p + \frac{3}{2}\sum_{p\le k} \frac{\lam_f(p)^2}{p} 
= \frac{13}{6} \sum_{p\le k} \frac 1p + \frac 32 \sum_{p\le k} \frac{\lam_f(p^2)}{p}.
$$
Using Lemma 2, and that $\sum_{p\le k}1/p = \log \log k +O(1)$, the Lemma follows.

\enddemo

\demo{Remark}  The {\sl ad hoc} inequality $2|x|\le \frac 23+\frac 32 x^2$ used above
was chosen for the simplicity of the statement in Lemma 3. 
If we write $L(1,\text{sym}^2 f) =(\log k)^{\theta}$ (so that $\theta \ge -1+o(1)$ by (2.2)), 
then a more complicated, but more precise, bound is 
$$ 
M_k(f) \ll (\log k)^{-2-\theta +2\sqrt{1+\theta+\epsilon}}.
$$
This follows upon using $2|x|\le \sqrt{1+\theta+\epsilon} +x^2/\sqrt{1+\theta+\epsilon}$ in 
the argument above.
\enddemo 

 \head 3.  Proof of Theorem 1 \endhead
 
 \noindent {\bf Case (i): The inner product with Maass cusp forms.}  Suppose that 
 $L(1,\text{sym}^2 f) \ge (\log k)^{-\frac {7}{15}}$.  Then Theorem 3 (i) gives 
 that $|\langle \phi F_k,F_k \rangle| \ll_{\phi} (\log k)^{-\frac 1{30}+\epsilon} $.
 
 Now suppose that $L(1,\text{sym}^2 f) < (\log k)^{-\frac {7}{15}}$.  Then Lemma 3 gives that 
 $M_k(f) \ll (\log k)^{-\frac{1}{15} +\epsilon}$, and using this in Theorem 2 (i), we obtain again 
 that $|\langle \phi F_k,F_k \rangle| \ll_{\phi} (\log k)^{-\frac 1{30}+\epsilon} $.

 In either case, the bound stated in Theorem 1 (i) holds. 
   \medskip
 \noindent {\bf Case (ii): The inner product with incomplete Eisenstein series.} We 
 begin by showing how Theorem 3 (ii) applies to incomplete Eisenstein series.  
 By Mellin inversion we 
 may write 
 $$ 
 E(z\mid \psi) = \frac{1}{2\pi i} \int_{\sigma-i\infty}^{\sigma+i\infty} {\hat \psi}(-s) E(z,s) ds, \tag{3.1}
 $$ 
 where $\sigma>1$ so that we are in the range of absolute convergence 
 of the Eisenstein series $E(z,s)$, and ${\hat \psi}$ denotes the Mellin transform 
 $$
 {\hat \psi}(s) = \int_0^{\infty} \psi(y) y^s \frac{dy}{y}.
 $$
 Since $\psi$ is smooth and compactly supported inside $(0,\infty)$ we 
 have that ${\hat \psi}$ is an analytic function, and repeated integration by 
 parts shows that $|{\hat \psi}(s)| \ll_{A,\psi} (1+|s|)^{-A}$ for any positive 
 integer $A$.  
We now move the line of integration in (3.1) to the line Re$(s)=\frac 12$.  
The pole of $E(z,s)$ at $s=1$ leaves a residue $\frac 3\pi$ and so 
$$ 
E(z \mid \psi) = \frac{3}{\pi} {\hat \psi}(-1) + \frac{1}{2\pi i} \int_{\frac 12-i\infty}^{\frac 12+i\infty} 
{\hat \psi}(-s) E(z, s) ds. \tag{3.2}
$$

From (3.2) it follows that 
 $$ 
 \langle E(\cdot \mid \psi ) F_k, F_k\rangle = \frac 3\pi {\hat \psi}(-1) + 
 \frac{1}{2\pi i} \int_{\frac 12-i\infty}^{\frac 12+i\infty} {\hat \psi}(-s) \langle  F_k E(\cdot, s), F_k \rangle \  ds.
 \tag{3.3}
 $$
 By unfolding we see that 
 $$ 
 \langle E(\cdot, \psi),1 \rangle = \int_{0}^{\infty} \int_{-\frac 12}^{\frac 12} \psi(y) \frac{dx \ dy}{y^2} 
 = {\hat \psi}(-1),
 $$ 
 and by Theorem 3 (ii) it follows that 
 $$ 
  \align
  \frac{1}{2\pi i} \int_{\frac 12-i\infty}^{\frac 12+i\infty} {\hat \psi}(-s) \langle   E(\cdot, s)F_k, F_k \rangle \  ds 
  &\ll \int_{-\infty}^{\infty} |{\hat \psi}(-\tfrac 12-it)| 
  \frac{(1+|t|)^2}{(\log k)^{1-\epsilon} L(1,\text{sym}^2 f)} dt\\
  &\ll \frac{(\log k)^{-1+\epsilon}}{L(1,\text{sym}^2 f)}.\\
  \endalign
  $$
  Using these in (3.3) we conclude that 
  $$ 
  \Big| \langle E(\cdot \mid \psi)F_k, F_k \rangle - \frac 3\pi \langle E(\cdot \mid \psi), 1 \rangle 
  \Big| \ll \frac{(\log k)^{-1+\epsilon}}{L(1,\text{sym}^2 f)}.
  $$
  Thus if  $L(1,\text{sym}^2 f) \ge (\log k)^{-\frac {13}{15}}$, 
  we obtain  the bound stated in Theorem 1 (ii).  
  
  Suppose now that $L(1,\text{sym}^2 f)\le (\log k)^{-\frac {13}{15}}$.  Then by Lemma 3 we have $M_k(f) \ll (\log k)^{-\frac {4}{15}+\epsilon}$, and by Theorem 2 (ii) (using Lemma 1 to estimate $R_k(f)$ there) we obtain 
  $$
  \Big| \langle   E(\cdot \mid \psi)F_k, F_k \rangle - \frac 3\pi \langle E(\cdot \mid \psi), 1 \rangle 
  \Big| \ll (\log k)^{\epsilon} M_k(f)^{\frac 12} \ll (\log k)^{-\frac 2{15} +\epsilon}, 
  $$
  so that the bound stated in Theorem 1 (ii) follows in this case also.

\enddocument